\documentclass[12pt,twoside]{article}
\usepackage{amsfonts}
\usepackage{theorem}
\usepackage{amsmath}

\newtheorem{theorem}{Theorem}[section]

\newtheorem{lemma}[theorem]{Lemma}

\newtheorem{remark}[theorem]{Remark}
\newtheorem{proposition}[theorem]{Proposition}

\newcommand{\openbox}{$\begin{array}{c}
\hspace*{-0.55em}\sqcap \hspace*{-0.60em}\\[-0.4em] \hline
\multicolumn{1}{c}{\hspace*{-0.60em}}\\[-0.8em]
\end{array}$}

\usepackage{enumerate}
\usepackage{graphics}
\usepackage{graphicx}
\usepackage{tikz}
\makeindex

\begin{document}

\centerline{\bf Remarks on Graphons}

\bigskip
\centerline{\bf Attila Nagy\footnote{This work was supported by the National Research, Development and Innovation Office – NKFIH, 115288.}}

\bigskip
\medskip

\noindent
\centerline{Department of Algebra}
\centerline{Budapest University of Technology and Economics}
\centerline{1521 Budapest, Pf. 91, Hungary}
\centerline{e-mail: nagyat@math.bme.hu}

\bigskip

\begin{abstract}
L. Lov\'asz and B. Szegedy proved in 2006 that
the limits of convergent graph sequences can be described by measurable symmetric functions $W: [0, 1]\times [0, 1]\to [0, 1]$ called graphons.
In our present paper we investigate the structure of the set of all graphons within the semigroup
$(\mathfrak{F}([0, 1]^2); \circ)$
of all fuzzy subsets of the unit square $[0,1]^2=[0, 1]\times [0, 1]$, where the operation $\circ$ is defined by: for every $f, g\in \mathfrak{F}([0,1]^2)$ and every $s\in [0,1]^2$,
$(f\circ g)(s)=\vee _{x\in [0,1]^2}(f(x)\wedge g(s))$.
\end{abstract}
 
\bigskip


 
 \bigskip

\section{Introduction and motivation}
Let $G_n$ be a sequence of finite simple graphs whose number of nodes tends to infinity. For every fixed finite simple graph $F$, let $hom(F,G_n)$ denote the number of all homomorphisms from $F$ into $G_n$, that is, the edge-preserving functions from $V(F)$ into $V(G_n)$. Put \[t(F,G_n)=\frac{hom(F,G_n)}{|V(G_n)|^{|V(F)|}}.\] Clearly, $t(F,G_n)$ is the probability that a random mapping from $V(F)$ into $V(G_n)$ should be a homomorphism. The sequence $G_n$ is called convergent if $lim_{n\to \infty}t(F,G_n)$ exists for every finite simple graph $F$. Let \[t(F)=lim_{n\to \infty}t(F,G_n).\] Then $t$ is a graph parameter, that is, a function on simple graphs that is invariant under isomorphism. In \cite{Lovasz1}, the authors given characterizations of graph parameters that arise in this manner; that is, the authors characterize  the set $\mathfrak{T}$ of graph parameters $t$ for which there is a convergent sequence of simple graphs $G_n$ such that $t(F)=lim_{n\to \infty}t(F,G_n)$ for every simple graph $F$.
In the characterization of $\mathfrak{T}$, the symmetric and measurable functions $W: [0,1]^2=[0, 1]\times [0, 1]\mapsto [0,1]$ called graphons play an important role. Recall that a function $W:[0,1]^2\mapsto [0,1]$ is said to be symmetric if $W(x,y)=W(y,x)$ is satisfied for all $x, y\in [0,1]$.
A graph is said to be {\it $k$-labelled} ($k$ is a positive integer) if the graph has $k$ nodes labelled by $1, 2, \dots , k$. For a $k$-labelled simple graph $F$ and a graphon $W$, the integral
\[t(F,W)=\int _{[0, 1]^k}\prod _{ij\in E(F)}W(x_i,x_j)dx_1dx_2\cdots dx_k\] is called the \textit{density of the graph $F$ in the graphon $W$} (\cite{Lovasz11}), where $E(F)$ denotes the set of all edges of $F$.
In \cite[Theorem 2.2]{Lovasz1} it was shown that a graph parameter $t$ belongs to $\mathfrak{T}$ if and only if there is a graphon $W$ such that $t(F)=t(F,W)$ for all simple graphs $F$.

A function of a non-empty set $S$ into the real unit interval $[0, 1]$ is called a {\it fuzzy subset} of $S$ (see \cite{Zadeh}).
By \cite{Wang} and \cite{Mordeson}, if $*$ is an associative operation on a non-empty set $S$, then
the set $\mathfrak{F}(S)$ of all fuzzy subsets of $S$ form a semigroup under the operation $\circ$ defined by the following way: for arbitrary $f, g\in \mathfrak{F}(S)$ and $s\in S$,
\begin{equation}\label{circ0}
(f\circ g)(s)=\begin{cases}
\vee _{s=x*y}(f(x)\wedge g(y)), & \text{if $s\in S^2$}\\
0, & \text{otherwise.}\end{cases}
\end{equation}
As every graphon is a fuzzy subset of the unit square $[0, 1]^2$, the following problem seems interesting from a semigroup theory perspective.

\medskip

\noindent
{\bf Problem}:
{\it If an associative operation $*$ is given on the unit square $[0, 1]^2$, what can we say about the structure of the set ${\cal W}_0$ of all graphons in the semigroup $(\mathfrak{F}([0, 1]^2); \circ )$? Is it true that
${\cal W}_0$ forms a substructure of $(\mathfrak{F}([0, 1]^2); \circ )$? If so, what kind of substructure is it?}

\medskip

In this paper we deal with this problem in a special case: the given associative operation $*$ on $[0, 1]^2$ satisfies the identity $(x, y)*(u, v)=(u, v)$.
A semigroup $(S; *)$ is called a {\it right zero semigroup} if it satisfies the identity $a*b=b$.
With this terminology, the above problem is examined in that case when $[0, 1]^2$ is a right zero semigroup.

We note that if $S$ is a non-empty set (and so it is a right zero semigroup), then the operation $\circ$ defined in (\ref{circ0}) has the following form:
\begin{equation}\label{circle}
(f\circ g)(s)=\vee _{x\in S}(f(x)\wedge g(s)).
\end{equation}

Throughout the paper, for a non-empty set $S$, $(\mathfrak{F}(S); \circ )$ will denote the semigroup in which the operation $\circ$ is defined by (\ref{circle}). Thus the purpose of this paper is to examine the structure of the set ${\cal W}_0$ of all graphons in the semigroup $(\mathfrak{F}([0, 1]^2); \circ )$.
Our studies consist of two parts. In Section~\ref{1:1} we describe the structure of the semigroup $(\mathfrak{F}(S); \circ )$ for an arbitrary non-empty set $S$, in Section~\ref{2:2} we focus on the semigroup $(\mathfrak{F}([0, 1]^2);\circ )$ and its subset ${\cal W}_0$.
A semigroup $S$ is called a {\it band} if every element $e$ of $S$ is an idempotent element, that is, $e^2=e$. A band satisfying the identity $axa=xa$ is called a {\it right regular band} (\cite{Petrich}).
In Section~\ref{1:1} we prove that if $S$ is an arbitrary non-empty set, then the semigroup $(\mathfrak{F}(S); \circ )$ is a right regular band (Theorem~\ref{rr}). In Section~\ref{2:2}, applying the above result for the right regular band $(\mathfrak{F}([0, 1]^2);\circ )$, we show that the set
${\cal W}_0$ of all graphons is a left ideal of $(\mathfrak{F}([0, 1]^2);\circ )$. By this result, if $W$ is a graphon and $f$ is a fuzzy subset of $[0, 1]^2$, then $f\circ W$ is a graphon.
Thus, for arbitrary simple graphs $F$, we can consider the densities $t(F;W)$ and $t(F; f\circ W)$ of $F$ in $W$ and in $f\circ W$, respectively.
In Section~\ref{2:2} we give an upper bound to $|t(F; W)-t(F; f\circ W)|$. In Theorem~\ref{upperbound} we show that
$|t(F; W)-t(F; f\circ W)|\leq |E(F)|(\sup(W)-\sup(f))\Delta (\{W> \sup(f)\})$, where $\Delta (\{W> \sup(f)\})$ denotes the area of the set
$\{W> \sup(f)\}=\{ (x, y)\in [0, 1]^2:\ W(x, y)>\sup(f)\}$.

\medskip

For notations and notions not defined here, we refer to the paper \cite{Lovasz1} and the books \cite{Clifford1}, \cite{Lovasz2}, \cite{Nagy}, and \cite{Petrich}.

\section{On the semigroup $(\mathfrak{F}(S); \circ )$, where $S$ is an arbitrary non-empty set}\label{1:1}
For a fuzzy subset $f$ and a subset $X$ of a non-empty set $S$, let $\sup _X(f)=\vee _{x\in X}f(x)$. Especially, let $\sup (f)=\sup _S(f)$. If $f$ and $g$ are arbitrary fuzzy subsets of $S$, then
let $g_{f}$ and $g^*_f$ denote the following fuzzy subsets of $S$: for an arbitrary $s\in S$, let
\[g_f(s)=\begin{cases}
\sup(f), & \text{if $g(s)> \sup(f)$}\\
g(s), & \text{otherwise}\end{cases}\]
and
\[g^*_f(s)=\begin{cases}
g(s)-\sup(f), & \text{if $g(s)> \sup(f)$}\\
0, & \text{otherwise.}\end{cases}\]

\begin{remark}\label{osszeg} \rm By the above definitions, $g_f+g^*_f=g$ for every fuzzy subsets $f$ and $g$ of a non-empty set $S$.
\end{remark}

\begin{remark}\label{gf=g} \rm Let $f$ and $g$ be arbitrary fuzzy subsets of a non-empty set $S$.
It is clear that $\sup(g)\leq \sup(f)$ implies $g(s)\leq \sup(f)$ for every $s\in S$ and so $g_f=g$.
In case $\sup(g)>\sup(f)$, there is an element $s\in S$ such that $g(s)>\sup(f)$ and so $g_f(s)=\sup(f)<g(s)$. Hence $g_f\neq g$.
Thus, for every fuzzy subsets $f$ and $g$ of $S$,  the equation $g_f=g$ holds if and only if $\sup(g)\leq \sup(f)$.
\end{remark}

\medskip

By Remark~\ref{gf=g}, the following lemma holds.

\begin{lemma}\label{sup} For arbitrary fuzzy subsets $f$ and $g$ of a non-empty set $S$, the equations $g_f=g$ and $f_g=f$ together hold if and only if $\sup(g)=\sup(f)$.
\end{lemma}

\medskip
The next lemma will be used in Lemma~\ref{measurable}.

\begin{lemma}\label{supp} If $f$ and $g$ are fuzzy subsets of a non-empty set $S$ such that $\sup (f)\leq \sup (g)$ then $\sup (g_f)= \sup (f)$ and $\sup (g^*_f)=\sup (g)-\sup (f)$.
\end{lemma}

\noindent
{\bf Proof}. By the definition of $g_f$ and $g^*_f$, it is obvious.\hfill\openbox

\begin{theorem} \label{fog} Let $S$ be a non-empty set. For every fuzzy subsets $f$ and $g$ of $S$, we have $f\circ g=g_f$.
\end{theorem}

\noindent {\bf Proof}.
Let $f$ and $g$ be arbitrary fuzzy subsets of a non-empty set $S$. By the above, $(\mathfrak{F}(S); \circ )$ is a semigroup. Let $s$ be an arbitrary element of $S$.
If $g(s)> \sup(f)$, then $f(x)\wedge g(s)=f(x)$ for every $x\in S$, and so
$(f\circ g)(s)=\vee _{x\in S}f(x)=\sup(f)$.
If $g(s)\leq \sup(f)$, then we have two subcases.

\noindent
Case 1: If $g(s)=\sup(f)$, then $f(x)\wedge g(s)=f(x)$ for all $x\in S$, and so
$(f\circ g)(s)=\vee _{x\in S}f(x)=\sup(f)=g(s)$.

\noindent
Case 2: If $g(s)<\sup(f)$, then there is an $x_0\in S$ such that $f(x_0)>g(s)$ and so $f(x_0)\wedge g(s)=g(s)$. Moreover, for arbitrary $x\in S\setminus \{ x_0\}$, we have
\[f(x)\wedge g(s)=\begin{cases}
g(s) ,& \textit{if $g(s)<f(x)$}\\
f(x) , & \textit{if $f(x)\leq g(s)$,}\end{cases}\] and so
$(f\circ g)(s)=(f(x_0)\wedge g(s))\vee(\vee _{x\in S\setminus \{ x_0\}}(f(x)\wedge g(s))=g(s)$. Summarizing our results, we get
\[(f\circ g)(s)=\begin{cases}
\sup(f), & \text{if $g(s)> \sup(f)$}\\
g(s), & \text{otherwise,}\end{cases}\] that is,
$(f\circ g)(s)=g_f(s)$, which proves our assertion.
\hfill\openbox

\medskip
A commutative band is called a \textit{semilattice}.
A congruence $\alpha$ on a semigroup $A$ is said to be a \textit{semilattice congruence} if the factor semigroup $A/\alpha$ is a semilattice. A semigroup $A$ is said to be \textit{semilattice indecomposable} if the universal relation is the only semilattice congruence on $A$. It is known (\cite{Tamura}) that every semigroup has a least semilattice congruence $\eta$; the classes of $\eta$ are semilattice indecomposable. By \cite[II.3.12. Proposition]{Petrich}, a band is a right regular band if and only if its $\eta$-classes are right zero semigroups.

\medskip

\begin{theorem} \label{rr} For an arbitrary non-empty set $S$,
the semigroup $(\mathfrak{F}(S); \circ )$ is a right regular band. The $\eta$-classes of $\mathfrak{F}(S)$ are right zero semigroups. Two fuzzy subsets $f$ and $g$ of $S$ are in the same $\eta$-class if and only if $\sup(f)=\sup(g)$.
\end{theorem}
\medskip

\noindent
{\bf Proof}.
Let $S$ be an arbitrary non-empty set. Then $S$ is a right zero semigroup, and so $(\mathfrak{F}(S); \circ )$ is a semigroup under the operation $\circ$ defined in (\ref{circle}), that is, $(f\circ g)(s)=\vee _{x\in S}(f(x)\wedge g(s))$ for every fuzzy subsets $f$ and $g$ of $S$ and every element $s\in S$. By Theorem~\ref{fog}, it is clear that $f\circ f=f$ for every $f\in \mathfrak{F}(S)$, and so $(\mathfrak{F}(S); \circ )$ is a band. Using also Theorem~\ref{fog}, we have $g\circ f\circ g=g\circ g_f$. As $\sup(g)\geq \sup(g_f)$, we have $g\circ g_f=g_f$.
Thus $g\circ f\circ g=g_f=f\circ g$. Hence $(\mathfrak{F}(S); \circ )$ is a right regular band.
Let $\eta$ denote the least semilattice congruence on $(\mathfrak{F}(S); \circ )$.
The $\eta$-classes of $(\mathfrak{F}(S); \circ )$ are right zero semigroups by \cite[II.3.12. Proposition]{Petrich}.
Let $f$ and $g$ be arbitrary fuzzy subsets of $S$. By \cite[II.1.1. Proposition]{Petrich}, $(f, g)\in \eta$ if and only if $f\circ g\circ f=f$ and $g\circ f\circ g=g$.
As $(\mathfrak{F}(S); \circ )$ is a right regular band, we have
$f\circ g\circ f=g\circ f$ and $g\circ f\circ g=f\circ g$. Thus
$(f, g)\in \eta$ if and only if $g\circ f=f$ and $f\circ g=g$. Using Theorem~\ref{fog},
$(f, g)\in \eta$ if and only if $f_g=f$ and $g_f=g$. By Lemma~\ref{sup}, we get
$(f, g)\in \eta$ if and only if $\sup(f)=\sup(g)$.
\hfill\openbox

\section{On the structure of the set of all graphons in the semigroup $(\mathfrak{F}([0, 1]^2); \circ )$}\label{2:2}

\medskip

Let $(S, {\cal A}, \mu )$ be a measurable space (\cite{Cohn}). For a fuzzy subset $h$ of $S$ and a real number $A$, let
$\{ h>A\}=\{ s\in S:\ h(s)>A\}$.
A fuzzy subset $h$ of $S$ is said to be \textit{measurable} if, for every real number $A$, the subset $\{ h>A\}$ of $S$ is measurable (that is, $\{ h>A\}\in {\cal A}$).

\begin{lemma} \label{measurable} Let $(S, {\cal A}, \mu )$ be a measurable space. Then, for an arbitrary fuzzy subset $f$ and an arbitrary measurable fuzzy subset $g$ of $S$, the fuzzy subsets $g_f$ and $g^*_f$ are measurable.
\end{lemma}
\medskip

\noindent
{\bf Proof}.
Let $f$ and $g$ be arbitrary fuzzy subsets of $S$ such that $g$ is measurable. If $\sup (f)\geq \sup (g)$, then $g_f=f\circ g=g$ and $g^*_f=0$. In this case the fuzzy subsets $g_f$ and $g^*_f$ are measurable.
Consider the case when $\sup (f)< \sup (g)$. Then $\sup (g_f)=\sup (f)$ and $\sup (g^*_f)=\sup (g)-\sup (f)$ by Lemma~\ref{supp}.
Let $A$ be an arbitrary real number.
It is easy to see that
\[\{g_f>A\}=\begin{cases}
\emptyset, & \text{if $A\geq \sup(f)$}\\
\{g>A\}, & \text{otherwise}\end{cases}\]
and
\[\{g^*_f>A\}=\begin{cases}
\emptyset, & \text{if $A\geq \sup(g)-\sup(f)$}\\
\{g>A+\sup(f)\}, & \text{if $0\leq A < \sup(g)-\sup (f)$}\\
S, & \text {if $A< 0$}\end{cases}\]
from which it follows that $g_f$ and $g^*_f$ are measurable fuzzy subsets of $S$.
\hfill\openbox

\medskip

A fuzzy subset $f$ of $[0, 1]^2$ is said to be \textit{symmetric} if $f(x,y)=f(y,x)$ is satisfied for all $x, y\in [0, 1]$.

\begin{lemma}\label{symleftid} If $f$ is an arbitrary fuzzy subset and $g$ is a symmetric fuzzy subset of $[0, 1]^2$, then
$g_f$ and $g^*_f$ are symmetric fuzzy subsets of $[0, 1]^2$.
\end{lemma}

\noindent
{\bf Proof}.
It is obvious by the definition of $g_f$ and $g^*_f$.
\hfill\openbox

\begin{lemma}\label{mindketto} If $W$ is a graphon and $f$ is a fuzzy subset of $[0, 1]^2$, then $W_f$ and $W^*_f$ are graphons.
\end{lemma}

\noindent
{\bf Proof}.
By Lemma~\ref{measurable} and Lemma~\ref{symleftid}, it is obvious.
\hfill\openbox

\medskip
The following theorem provides an answer to the question raised in Problem in the case, where the given operation $\cdot$ on $[0, 1]^2$ satisfies the identity $a\cdot b=b$.

\begin{theorem}\label{leftideal} The set ${\cal W}_0$ of all graphons is a left ideal of the right regular band $(\mathfrak{F}([0, 1]^2); \circ )$ of all fuzzy subsets of $[0, 1]^2$. Thus the semigroup $({\cal W}_0; \circ )$ of all graphons is a right regular band, and so it is a semilattice $I$ of right zero subsemigroups $S_i$ ($i\in I$). Two graphons $W_1$ and $W_2$ are in the same $S_i$ if and only if $\sup(W_1)=\sup(W_2)$.
\end{theorem}

\noindent
{\bf Proof}.
Let $W$ be a graphon and $f$ be a fuzzy subset of $[0, 1]^2$. By Theorem~\ref{fog}, $f\circ W =W_f$. Then $f\circ W$ is a graphon by Lemma~\ref{mindketto}. Thus the set ${\cal W}_0$ of all graphons is a left ideal of the semigroup $(\mathfrak{F}([0, 1]^2); \circ )$ of all fuzzy subsets of $[0, 1]^2$. By Theorem~\ref{rr}, the semigroup $(\mathfrak{F}([0, 1]^2); \circ )$ and so its subsemigroup $({\cal W}_0; \circ )$ is a right regular band. Moreover, the $\eta$-classes of ${\cal W}_0$ are right zero semigroups; two graphons $W_1$ and $W_2$ are in the same $\eta$-class if and only if $\sup(W_1)=\sup(W_2)$.
\hfill\openbox

\medskip

Let $\sigma$ denote the equivalence relation on the set ${\cal W}_0$ of all graphons defined by
$(W_1, W_2)\in \sigma$ if and only if $W_1=W_2$ almost everywhere in $[0, 1]^2$.

\begin{proposition} The equivalence relation $\sigma \cap \eta$ is a congruence on the right regular band $({\cal W}_0; \circ )$ of all graphons, where $\eta$ is the least semilattice congruence on $({\cal W}_0; \circ )$.
\end{proposition}

\noindent
{\bf Proof}.
Let $W_1$ and $W_2$ be two graphons with $(W_1, W_2)\in \sigma \cap \eta$. Then, using Theorem~\ref{leftideal}, we have
$\sup(W_1)=\sup(W_2)$ and $W_1=W_2$ almost everywhere in $[0, 1]^2$. Let $W$ be an arbitrary graphon. As
$\sup(W_1)=\sup(W_2)$, we have $W_1\circ W=W_2\circ W$. Thus $(W_1 \circ W, W_2\circ W)\in \sigma \cap \eta$. Hence $\sigma \cap \eta$ is a right congruence on $({\cal W}_0; \circ )$.
Let $T=\{(x, y)\in [0, 1]^2|\ W_1(x, y)\neq W_2(x, y)\}$. As $(W_1, W_2)\in \sigma$, the area of $T$ is $0$. It is clear that $\{ (x, y)\in [0, 1]^2: (W\circ W_1)(x, y)\neq (W\circ W_2)(x, y)\} \subseteq T$
and so $(W\circ W_1, W\circ W_2)\in \sigma$. As $(W_1, W_2)\in \eta$ and $\eta$ is a congruence on $({\cal W}_0;\circ )$, we have
$(W\circ W_1, W\circ W_2)\in \eta$. Thus $(W\circ W_1, W\circ W_2)\in \sigma \cap \eta$ and so $\sigma \cap \eta$ is a left congruence on $({\cal W}_0; \circ )$. Thus
$\sigma \cap \eta$ is a congruence on $({\cal W}_0; \circ )$.
\hfill\openbox

\medskip

Let $W$ be a graphon and $f$ a fuzzy subset of $[0, 1]^2$. By Theorem~\ref{leftideal}, $f\circ W$ is a graphon.
Thus, for arbitrary simple graphs $F$, we can consider the densities $t(F;W)$ and $t(F; f\circ W)$ of $F$ in $W$ and $f\circ W$, respectively. The next theorem gives an upper bound to $|t(F; W)-t(F; f\circ W)|$.

\begin{theorem}\label{upperbound} Let $W$ be an arbitrary graphon. Then, for an arbitrary fuzzy subset $f$ of $[0, 1]^2$ and an arbitrary finite simple graph $F$,
\[|t(F; W)-t(F; f\circ W)|\leq |E(F)|(\sup(W)-\sup(f))\Delta (\{W> \sup(f)\}),\] where $E(F)$ denotes the set of all edges of $F$ and
$\Delta (\{W> \sup(f)\})$ denotes the area of the set
$\{W> \sup(f)\}=\{ (x, y)\in [0, 1]^2:\ W(x, y)>\sup(f)\}$.
\end{theorem}

\noindent
{\bf Proof}.
Let $W$ be an arbitrary graphon and $f$ an arbitrary fuzzy subset of $[0, 1]^2$. By Theorem~\ref{leftideal}, $f\circ W$ is a graphon.
If $\sup(W)\leq \sup(f)$, then
$W=f\circ W$ and $\{W> \sup(f)\}=\emptyset$. Thus
$|t(F; W)-t(F; f\circ W)|=0=|E(F)|(\sup(W)-\sup(f))\Delta (\{W> \sup(f)\})$.
Consider the case when $\sup(W)>\sup(f)$. By Remark~\ref{osszeg}, $W-(f\circ W)=W^*_f$. As $W$ is a graphon, $W_f=f\circ W$ and $W^*_f$ are graphons by Lemma~\ref{mindketto}. Thus $W$, $f\circ W$ and $W^*_f$ are integrable functions on $[0, 1]^2$. Using \cite[Lemma 4.1]{Lovasz1}, $|t(F; W)-t(F; f\circ W)|\leq |E(F)|\cdot ||W^*_f||_0$,
where $||W^*_f||_0=\sup_{A\subseteq [0, 1]\atop B\subseteq [0, 1]}\left | \int _A\int_BW^*_f(x, y)dxdy\right |$. As $W^*_f$ is a non-negative function,
$||W^*_f||_0=||W^*_f||_1$, where
$||W^*_f||_1=\int _0^1\int _0^1|W^*_f(x, y)|dxdy$. Thus
$|t(F; W)-t(F; f\circ W)|\leq |E(F)|\cdot ||W^*_f||_1$.
As $W^*_f(x, y)=0$ for all $(x, y)\in [0, 1]^2\setminus \{ W>\sup(f)\}$,
we have
$||W^*_f||_1=\int_0^1\int_0^1W^*_f(x, y)dxdy\leq (\sup(W)-\sup(f))\Delta (\{W>\sup(f)\}$,
because $\sup(W^*_f)= \sup(W)-\sup(f)$ by Lemma~\ref{supp}.
Consequently
$|t(F; W)-t(F; f\circ W)|\leq |E(F)|(\sup(W)-\sup(f))\Delta (\{W> \sup(f)\})$.
\hfill\openbox

\end{document}